\documentclass[10pt,twoside,reqno]{amsart}
\usepackage{amsthm}
\usepackage{amsmath}
\usepackage{amssymb}
\usepackage{amsfonts}
\usepackage{amscd}
\usepackage{amsxtra}

\usepackage{eucal}
\usepackage{mathrsfs}
\usepackage{nicefrac}
\usepackage{bbold}

\usepackage[ruled]{algorithm} 
\usepackage{algorithmic}

\newtheorem{theorem}{Theorem}[section]

\newtheorem*{proposition*}{Proposition}

\newtheorem*{theorem*}{Theorem}
\newtheorem*{lemma*}{Lemma}
\newtheorem*{conjecture*}{Conjecture}

\theoremstyle{remark}

\newtheorem*{remark*}{Remark}
\newtheorem{remsTh}[theorem]{Remarks}

\newcommand{\enumTi}[1]{\renewcommand{\theenumi}{#1}}

\renewcommand{\emph}{\sl}
\newcommand{\stepnr}[1]{{\footnotesize #1}}

\newcommand{\lt}{\left}
\newcommand{\rt}{\right}

\newcommand{\abs}[1]{{\lt\lvert{#1}\rt\rvert}}
\newcommand{\sabs}[1]{{\lvert{#1}\rvert}}

\newcommand{\nfrac}[2]{{\nicefrac{#1}{#2}}}

\newcommand{\cplmt}{\complement}

\newcommand{\ZZ}{\mathbb{Z}}
\newcommand{\QQ}{\mathbb{Q}}
\newcommand{\RR}{\mathbb{R}}

\newcommand{\AlgoFloatname}{Algorithm}
\floatname{algorithm}{\AlgoFloatname}

\renewcommand{\emph}[1]{{\it #1}}

\begin{document}

\title{Odd minimum cut sets and $b$-matchings revisited}%

\author{Adam N. Letchford}%
\address{Adam N. Letchford\\Department of Management Science\\Lancaster University\\Lancaster LA1 4YW\\England.}
\email{A.N.Letchford@lancaster.ac.uk}

\author{Dirk Oliver Theis${}^*$}
\address{Dirk Oliver Theis\\Faculty of Mathematics and Computer Science\\University of Heidelberg\\Germany.}
\email{theis@uni-heidelberg.de}
\thanks{${}^*$Supported within project RE~776/9-1 of the Deutsche Forschungsgemeinschaft (DFG)}

\begin{abstract}
The famous Padberg--Rao separation algorithm for $b$-matching polyhedra can be implemented to run in
${\mathcal O}(|V|^2|E| \log (|V|^2/|E|))$ time in the uncapacitated case, and in
${\mathcal O}(|V||E|^2 \log (|V|^2/|E|))$ time in the capacitated case. We give a new and simple algorithm
for the capacitated case which can be implemented to run in ${\mathcal O}(|V|^2|E| \log (|V|^2/|E|))$ time.
\end{abstract}

\maketitle


\section{Introduction}

Let $G=(V,E)$ be an undirected graph, let $b \in \ZZ_+^V$ be a vector of vertex capacities and let
$u \in \ZZ_+^E$ be a vector of edge capacities. A \emph{$u$-capacitated $b$-matching} is a family of
edges, possibly containing multiple copies, such that:
\begin{itemize}
\item for each $i \in V$, there are at most $b_i$ edges in the family incident on $i$;
\item at most $u_e$ copies of edge $e$ are used.
\end{itemize}
If we define for each edge $e$ the integer variable $x_e$, representing the number of times $e$
appears in the matching, then the incidence vectors of $u$-capacitated $b$-matchings are the
solutions to:
\begin{align}    
  {\textstyle\sum\nolimits_{e\in\delta(i)}}\ x_e \leq b_i, &\quad\text{for all }i \in V &   \label{eq:degree}\\
  0 \leq x_e \leq u_e,   &\quad\text{for all }e \in E &  \label{eq:bounds}\\
  x_e \in \ZZ,           &\quad\text{for all }e \in E.&  \label{eq:int}
\end{align}
Here, as usual, $\delta(i)$ represents the set of vertices incident on $i$.

The convex hull in $\RR^E$ of solutions to (\ref{eq:degree}) - (\ref{eq:int}) is called the
\emph{$u$-capacitated $b$-matching polytope}.  Edmonds and Pulleyblank (see \cite{Edmonds65} and
\cite{PulleyblankPhD}) gave a complete linear description of this polytope. It is described by the
\emph{degree inequalities} (\ref{eq:degree}), the \emph{bounds} (\ref{eq:bounds}) and the following
{\em blossom inequalities}:
\begin{multline} \label{eq:blossom}
  \sum\nolimits_{e\in E(W)} x_e  + \sum\nolimits_{f\in F} x_f
  \leq \left\lfloor \frac{b(W) + \sum_{f\in F}u_f}{2} \right\rfloor,\\
  \text{for all }W \subset V\text{, }F \subset \delta(W)
  \text{ with }b(W) + {\textstyle \sum_{f\in F}}\ u_f \text{ odd}.
\end{multline}
Here, $E(W)$ (respectively, $\delta(W)$) represents the set of edges with both end-vertices
(respectively, exactly one end-vertex) in $W$, $b(W)$ denotes $\sum_{i \in W} b_i$.

An important special case is where the upper bounds $u_e$ are not present (or, equivalently, $u_{ij}
\geq \max \{b_i, b_j\}$ for all $\{i,j\} \in E$).  The associated (uncapacitated) $b$-matching
polytope is described by the degree inequalities, the non-negativity inequalities $x_e \geq 0$ for
all $e \in E$, and the {\em simplified blossom inequalities}
\begin{equation} \label{eq:simple}
\sum\nolimits_{e\in E(W)} x_e \leq \left\lfloor \frac{b(W)}{2} \right\rfloor,
\quad\text{ for all }W \subset V\text{ with }b(W) \text{ odd}.
\end{equation}

In their seminal paper, \cite{PadRao} devised a combinatorial, polynomial-time \emph{separation
algorithm} for $b$-matching polytopes.  A separation algorithm is a procedure which, given a
rational vector $x^* \in \QQ^E$ lying outside of the polytope, finds a linear inequality which is
valid for the polytope yet violated by $x^*$. Clearly, testing if a degree inequality or bound is
violated can be performed in linear time, so the main contribution of \cite{PadRao} is to identify
violated blossom inequalities.

For \emph{uncapacitated $b$-matching,} Padberg \& Rao reduce the separation problem to the
computation of a minimum $T$-cut, for which they give a generic algorithm, see
Algorithm~\ref{algo:PadraoMOC}.  We will give the definition of the minimum $T$-cut problem in the
next section.  Abbreviating $n:=\abs V$ and $m := \abs E$, this algorithm involves the solution of
up to $n-1$ maximum flow problems on a graph with $n+1$ vertices and $n+m$ edges. Using the
well-known \emph{pre-flow push} algorithm \cite{GoldbergTarjan} to solve the max-flow problems, this
leads to an overall running time of ${\mathcal O}(n^2m \log\nfrac{n^2}{m})$.

\begin{algorithm}[t]
  \caption{Minimum T-cut \cite{PadRao}}\label{algo:PadraoMOC}
  \flushleft
  \begin{algorithmic}[1]
    \REQUIRE~\\{ Graph $G$, set $T\subset V$, and weights $c \in\QQ_+^E$. }%
    \ENSURE~\\{ A minimum T-cut. }
    \\[1mm]
    \STATE{Compute a cut-tree for the graph $G$ with weights $c$ and terminal vertex set $T$. }
    \FOR{each of the $n-1$ edges of the cut-tree}{%
      \STATE Let $\delta(U)$ denote the cut induced by the cut-tree edge.
      \STATE{\textit{Check the cut:}\\
        Compute the parity $\sabs{T\cap U} \bmod 2$ and the weight $c(U)$ of the cut}.\\
      \STATE If adequate, store $U$.
    }\ENDFOR%
    \STATE Output the best $T$-cut $U$.
  \end{algorithmic}
\end{algorithm}%

The Padberg-Rao separation algorithm for \emph{capacitated $b$-matching,} however, is substantially
more time-consuming.  It involves the computation of a minimum $T$-cut on a special graph, the
so-called {\em split graph,} which has up to $n+m+1$ vertices and up to $2m+n$ edges.  Up to $n+m-1$
maximum flow problems may be required to be computed.  Using the pre-flow push algorithm, this leads
to a worst-case running time of ${\mathcal O}(m^3 \log n)$.
%
In 1987, \cite{GroetschelHolland87} observed that the above-mentioned max-flow problems can in fact
be carried out on graphs with only ${\mathcal O}(n)$ vertices and ${\mathcal O}(m)$ edges.  Although the
idea behind this is simple, it reduces the overall running time for the capacitated case to
${\mathcal O}(nm^2 \log\nfrac{n^2}{m})$.

In this paper, we propose a new separation algorithm for the capacitated case whose running time is
the same as that for the uncapacitated case. As well as being faster than the Padberg-Rao and
Gr\"otschel-Holland approaches, the new algorithm is much simpler and easier to implement. It also
has a surprisingly simple proof of correctness.

Our results also apply to the case of {\em perfect capacitated $b$-matchings.}

As well as being of interest in the context of matching, the algorithm has an important application
to the \emph{Traveling Salesman Problem} (TSP). The special blossom inequalities obtained when $b_i
= 2$ for all $i$ and $u_e = 1$ for all $e$ are valid for the TSP, and facet-inducing under mild
conditions, see \cite{GroePadb79a}, \cite{GroePadb79b}.  
Thus we obtain a faster exact separation algorithm for the TSP as a by-product.  In fact, the
algorithm is applicable to a general class of cutting planes for integer programs, called
$\{0,\nfrac12\}$-Chv\'atal-Gomory cuts, see \cite{CaprFisch96}.

Parts of the contents of this paper appeared in the proceedings of the Xth IPCO conference
\cite{LetchReineltTheis04a}.  However, the proof of correctness of the algorithm is now {\sl
  substantially} facilitated.

\begin{algorithm}[t]
  \caption{Blossom minimization}\label{algo:AdamBlossomMin}
  \flushleft
  \begin{algorithmic}[1]
    \REQUIRE~\\{ Graph $G$, set $T\subset V$, and weights $c,c' \in\QQ_+^E$. } %
    \ENSURE~\\{ A minimum blossom. }
    \\[1mm]
    \STATE{ Compute a cut-tree for $G$ with weights $\min(c,c')$ and
      terminal vertex set $V$. }
    \FOR{each of the $n-1$ edges of the cut-tree}{%
      \STATE Let $\delta(U)$ denote the cut induced by the cut-tree edge.
      \STATE\textit{Check the cut:}\\
      Compute $\beta(U)$ as in \eqref{eq:min-beta}.
      \STATE If adequate, store $U$ along with the arg-min $F$.
    }\ENDFOR%
    \STATE Output the best blossom $(U,F)$.
  \end{algorithmic}
\end{algorithm}


\section{Algorithms for minimum $T$-cut and blossom minimization}

Given a graph $G=(V,E)$, an even-cardinality set $T\subset V$ and non-negative rational
edge-capacities $c \in \QQ_+^E$, the {\em minimum $T$-cut problem} asks for an {\em odd cut}
$(U,\cplmt U)$ (where $\cplmt U$ is the complement of $U$ in the vertex set) such that the set
$U\subset V$ is {\em $T$-odd}, i.e., $\sabs{T\cap U}$ is an odd number, and which minimizes, subject
to this condition, the submodular function
\begin{equation*}
  U \mapsto c(U) := \sum_{e\in \delta(U)} c_e.
\end{equation*}
In 1982, Padberg \& Rao gave the first polynomial-time combinatorial algorithm for computing a
minimum T-cut, see Algorithm~\ref{algo:PadraoMOC}.  The key ingredient is the computation of a
Gomory-Hu cut-tree \cite{GomoryHu1961} in step~\stepnr1.
Given a graph $G=(V,E)$, a set $X\subset V$, and non-negative rational vector of edge-capacities
$c\in\QQ_+^E$, a {\em cut-tree with terminal vertex set $X$} for $G$ and $c$ consists of a mapping
$\pi\colon V\to X$ with $\pi(x)=x$ for all $x\in X$, and an adjacency relation $\sim$ on the set
$X$.  (We adopt the convention that the edges of $G$ will be denoted by $xy$, and the edges of the
cut-tree by $x\sim y$.)  The adjacency relation shall make the set of terminal vertices into a tree.
An additional condition is required to hold.  Deleting an edge $x\sim y$ of the cut-tree partitions
the set $X$ into two sets $X_x$ and $X_y$, and thus defines a cut $(U,\bar U)$ in $G$ by letting $U
:= \pi^{-1}(X_x)$ and $\bar U :=\pi^{-1}(X_y)$.  We call this the cut {\em induced} by the edge
$x\sim y$ of the cut-tree.  Now, the condition which is required is the following:
\begin{equation}  \label{eq:cut-tree-ax}
  \text{\parbox{11cm}{\flushleft%
      for $x,y\in X$ with $x\sim y$, the cut induced by this edge of the cut-tree shall be
      a minimum $(s,t)$-cut in $G$ with respect to the capacities $c$.
    }}
\end{equation}
With the algorithm given by Gomory \& Hu, a cut-tree can be computed in time
$O(\sabs{X}nm\log\nfrac{n^2}{m})$.

In Algorithm~\ref{algo:PadraoMOC}, the time for ``checking the cut'' in step~\stepnr4 is negligible
(the values $c(U)$ even come for free with the Gomory-Hu algorithm), and hence the Padberg-Rao
method for computing a minimum $T$-cut runs in time $O(\sabs{T}nm\log\nfrac{n^2}{m})$, as mentioned
in the introduction.

Now we come to the blossom separation algorithm of Padberg \& Rao \cite{PadRao}.  Reformulating and
generalizing, we say that a {\em blossom} is a pair $(U,F)$ consisting of a set of vertices
$U\subset V$ and a set of edges $F\subset \delta(U)$ with the property that $\sabs{T\cap U}+\sabs F$
is an odd number.  Then, if two non-negative rational weight vectors $c,c'\in\QQ_+^E$ are given for
the edges of $G$, the blossom separation problem is equivalent to the problem of producing a blossom
whose {\em value}
\begin{equation*}
  \beta(U,F) := \sum_{e\in \delta(U)\setminus F} c_e + \sum_{f\in F} c'_f
\end{equation*}
is strictly less than, one, if it exists.  For the sake of completeness, we describe how this
equivalence is established.  Padberg \& Rao \cite{PadRao} introduced, for each $u \in V$, the term
$s_u := b_u - \sum_{e\in\delta(i)} x_e$, which is the \emph{slack} of the corresponding degree
inequality computed with respect to a given vector $x$.  Then they showed that the blossom
inequality \eqref{eq:blossom} can be re-written in the form:
\begin{equation} \label{eq:oddcut}
  \sum_{u \in W} s_u
  + \sum\nolimits_{e\in\delta(i)} x_e
  +  \sum_{e \in F} (u_e - x_e)
  \geq 1.
\end{equation}
To decide if, for a given $x$, sets $W$ and $F$ exist which violate \eqref{eq:oddcut}, we define, in
a canonical and straight forward manner, a graph $G^*$, capacities $c$ and $c'$ and a set $T$ of
vertices of $G^*$, in such a way that a blossom with value strictly less than one gives rise to a
violated inequality \eqref{eq:oddcut} and vice-versa.  Let $G^*$ be constructed by adding a new
vertex $v$ to $G = (V,E)$ and connecting it with an edge $vu$ to every $u\in V$.  Then for each
$e\in E$, we let
\begin{equation*}
  (c_e,c'_e) :=
  \begin{cases}
    (x_e,u_e-x_e)              &\text{if $u_e$ is odd}\\
    (\min(x_e,u_e-x_e),\infty) &\text{if $u_e$ is even}
  \end{cases}
\end{equation*}
For the edges $vu$ of $G^*$, we let $c_{uv} := s_u$ and $c'_{uv} := \infty$.  Finally, we define $T$
as the set of all vertices $u$ for which the value $b_u$ is odd, and we let $v\in T$ iff $\sum_u b_u$
is odd.  Now it is easy to see that for each blossom $(U,F)$ in $G^*$ with $v\not\in U$, the
inequality \eqref{eq:oddcut} with $W:= U\cap V$ is violated by $1-\beta(U,F)$.  Note that
$\beta(U,F)=\beta(\cplmt U,F)$.

As mentioned above, the blossom separation Algorithm of Padberg \& Rao \cite{PadRao} is considerably
more complex than the minimum $T$-cut algorithm.  It requires to construct a special graph $\hat G$
with $m+n$ vertices and $2m$ edges, on which then a minimum $T$-cut is computed.

We now give an algorithm for what we call the {\em blossom minimization problem:} given $G$, $T$ and
$c,c'$ as above, find a blossom $(U,F)$ which minimizes $\beta(U,F)$.  The blossom minimization
algorithm is displayed as Algorithm~\ref{algo:AdamBlossomMin}.

For fixed $U\subset V$, it has been observed by Padberg \& Rinaldi \cite{PadRinaldiSepaTSP} that
\begin{equation}   \label{eq:min-beta}
  \beta(U) := 
  \min \Bigl\{  \beta(U,F) \Bigm\vert F\subset\delta(U)\text{, }\sabs{T\cap U}+\sabs F \text{ odd}   \Bigr\}
\end{equation}
can be computed in time $O(\sabs{\delta(U)})$ by first tentatively taking $F := \{e\in\delta(U)\mid
c'_e < c_e\}$.  Now if $\sabs{T\cap U}+\sabs{F}$ is odd, we have found a minimizing $F$.  Otherwise,
find $f\in\delta(U)$ minimizing $\lvert c_f-c'_f\bigr\rvert$ over $f\in\delta(U)$, because then the
symmetric difference of $F$ and $\{f\}$ minimizes $\beta(U,\cdot)$.

This implies that the loop \stepnr2--\stepnr6 in Algorithm~\ref{algo:AdamBlossomMin} runs in time
$O(n^2)$ and that the running time of Algorithm~\ref{algo:AdamBlossomMin} is dominated by the
computation of the cut-tree in step~\stepnr1, which amounts to $O(n^2m\log\nfrac{n^2}{m})$.

The similarity between the Padberg-Rao minimum $T$-cut Algorithm~\ref{algo:PadraoMOC} and our
blossom minimization Algorithm~\ref{algo:AdamBlossomMin} is striking.  Moreover, in the next
section, we give a short and elegant proof of correctness of Algorithm~\ref{algo:AdamBlossomMin},
which is similar to a proof of correctness of Algorithm~\ref{algo:PadraoMOC} given by Rizzi
\cite{Rizzi02}.  At this point, we might note that $\beta(\cdot)$, unlike $c(\cdot)$, is not in
general submodular.


\section{A simple proof of the correctness of Algorithm~\ref{algo:AdamBlossomMin}}

Let a cut-tree for $G$ with terminal vertex set $X$ be given, where $X\supset T$.  We say that an
edge $x\sim y$ of the cut-tree is {\em $T$-odd,} if the sets of the bipartition of $X$ defined by
$x\sim y$ are $T$-odd.  Thus, the set of $T$-odd edges of the cut-tree form what is called a
$T$-join, and an edge in the cut-tree induces a $T$-cut in $G$ if and only if the edge is $T$-odd.
The next theorem is the keystone of the correctness of Algorithm~\ref{algo:PadraoMOC}.  For the sake
of clarity, we repeat the proof of \cite{Rizzi02}.

\begin{theorem}[\cite{PadRao}]          \label{thm:padrao-rizzi}
  One of the $T$-odd edges of the cut-tree induces a minimum $T$-cut in $G$.
\end{theorem}

\begin{proof}
  Let $U$ be a minimum $T$-cut. Now $U$ is a $T$-odd set, hence there exists an odd number of
  $T$-odd cut-tree edges leaving $T\cap U$.  Let $x\sim y$ be one of them, and let $S$ be the
  minimum $(x,y)$-cut it induces by \eqref{eq:cut-tree-ax}.  Since $U$ is an $(x,y)$-cut, we have
  $c(S)\le c(U)$, and since $x\sim y$ is an $T$-odd edge, $S$ defines a minimum $T$-cut.
\end{proof}

Now we come to the proof of correctness of Algorithm~\ref{algo:AdamBlossomMin}.

\begin{theorem}
  One of the edges of of the cut-tree computed in Algorithm~\ref{algo:AdamBlossomMin} induces the a
  set $U$ which minimizes $\beta(\cdot)$.
\end{theorem}
\begin{proof}
  Let $U$ be a set which minimizes $\beta(\cdot)$.  Further, define the set $T'$ as the symmetric
  difference of $T$ with all sets $\{u,v\}$ for all $e=uv\in E$ $c'_e < c_e$.
  
  Case 1: $U$ is $T'$-odd.  The proof of Theorem~\ref{thm:padrao-rizzi} shows that there exists a
  $T'$-odd edge of the cut-tree which induces a minimizer of $\beta(\cdot)$.
  
  Case 2: $U$ is not $T'$-odd.  Let $f=x'y'\in\delta(U)$ have the minimal value of $\sabs{c_f-c'_f}$
  among all edges in $\delta(U)$.  On the path from $x'$ to $y'$ in the cut-tree, at least one edge
  $x\sim y$ has one end in $U$ and the other not in $U$.  Let $S$ be the minimum $(x,y)$-cut defined
  by this edge.  Abbreviating $w := \min(c,c')$, we then have
  \begin{equation*}
    \beta(U) = w(U) + \abs{c_f-c'_f} \ge w(S) + \abs{c_f-c'_f} \ge \beta(S).
  \end{equation*}
  The first inequality holds since $U$ is an $(x,y)$-cut.  As for the second, if $S$ is $T'$-odd,
  then $S$ minimizes $\beta$ since $w(S) \le w(S) + \sabs{c_f-c'_f} \le \beta(U)$; but if
  $\sabs{T'\cap S}$ is even, then $(S,\{f\})$ is a blossom whence $w(S) + \sabs{c_f-c'_f} =
  \beta(S,\{f\}) \ge \beta(S)$.
\end{proof}


\end{document}